\documentclass{article}
\usepackage{lscape, amsmath, graphicx, enumerate, pbox, hyperref}
\title{On \'Sankara Varman's (correct) and M\=adhava's (incorrect) values for the circumferences of circles}
\author{{\bf V. N. Krishnachandran}}
\allowdisplaybreaks
\date{}
\begin{document}
\maketitle

%\tableofcontents
%\newpage
%\thispagestyle{empty}
\begin{abstract}
This paper examines what computational procedures
\'Sankara Varman (1774-1839) and Sangamagrama M\=adhava (c. 1340 - 1425), astronomer-mathematicians of the Kerala school, might have used to arrive at their respective values for the circumferences of certain special circles (a circle of diameter $10^{17}$ by the former and a circle of diameter $9\times 10^{11}$ by the latter). It is shown that if we choose the M\=adhava-Gregory series for $\frac{\pi}{6}=\arctan \frac{1}{\sqrt{3}}$
   to compute $\pi$ and then use it compute the circumference of a circle of diameter $10^{17}$ and  perform the computations by ignoring the fractional parts in the results of every  operation we get the value stated by \'Sankara Varman. It is also shown 
   that, except in an unlikely case, none of the series representations of $\pi$ attributed to M\=adhava produce the value for the circumference attributed to him. The question  how M\=adhava did arrive at his value still remains unanswered.
\end{abstract}
\qquad\\[1mm]
\section{Introduction}
\'Sankara Varman (1774-1839), an astronomer-mathematician of the Kerala school of astronomy and mathematics, in his {\em Sadrtnam\=la} has given a  measure of the circumference of a circle of diameter $10^{17}$. This value is  correct to the true value of the circumference rounded to the nearest integer.
 Earlier mathematicians of the Kerala school like N\=ilaka\d{nt}ha Soma\=yj\=i and \'Sankara V\=ariar have quoted Sangamagrama M\=adhava (c. 1340 - 1425), the founder of the Kerala school, as stating that the circumference of a circle of diameter $9\times 10^{11}$ is $2827433388233$.  But the true value of the circumference of this circle when rounded to the nearest integer is  $2827433388231$ and so M\=adhava's value is in error by 2 units. This paper examines what computational procedures
these  astronomer-mathematicians might have used to arrive at their respective values for the circumferences of those special circles. It is shown that if we choose the M\=adhava-Gregory series for $\frac{\pi}{6}=\arctan \frac{1}{\sqrt{3}}$
   to compute $\pi$ and then use it compute the circumference of a circle of diameter $10^{17}$ and  perform the computations by ignoring the fractional parts in the results of every  operation we get the value stated by \'Sankara Varman. It is also shown 
   that, except in an unlikely case, none of the series representations of $\pi$ attributed to M\=adhava produce the value for the circumference attributed to him. The question  how M\=adhava did arrive at his value still remains unanswered.
   
In Section 2 we have a brief note on the life and work of \'Sankara Varman. In the next three sections we discuss \`sankara Varman's value. We discuss it first because the value is correct! In Section three we present the Sanskrit verses with their English translations given in {\em Sadratnam\=al\=a} which give the value of $\pi$ correct to 17 decimal places. Some issues like why a circle of diameter $10^{17}$, what rounding method, etc. connected with the computation of this value are discussed in Section 4. Sections 5 gives the complete details of the possible methods of computation all of which can be carried out by hand. In the last two sections we discuss the details of the computation of M\=adhava's value by using the various infinite series representations of the circumferences all attributed to M\=adhava. 
\section{\'Sankara Varman}
\'Sankara Varman (1774–1839), author of {\em Sadratnam\=ala} which is a  treatise on astronomy and mathematics composed in 1819, was an astronomer-mathematician belonging to the Kerala school of astronomy and mathematics. He is considered as the last notable figure in the long line of illustrious astronomers and mathematicians in the Kerala school of astronomy and mathematics beginning with M\=adhava of Sangamagr\=ama (c. 1340 - c. 1425). {\em Sadratnamala} was composed in the form of verses in Sanskrit in the traditional style followed by members of the Kerala school even though India had been introduced to the western style of writing books in mathematics at that time. 

\'Sankara Varman was a close acquaintance of C. M. Whish (1794 – 1833) a civil servant of East India Company attached to its Madras establishment. C. M. Whish was the first person  to bring to the attention of the European mathematical scholarship the mathematical accomplishments  of the Kerala school of astronomy and mathematics (see \cite{MacTutor2}).  Whish spoke of him and his work thus (see \cite{Whish}: ``The author of Sadratnamalah is SANCARA VARMA, the younger brother of the present Raja of Cadattanada near Tellicherry, a very intelligent man and acute mathematician. This work, which is a complete system of Hindu astronomy, is comprehended in two hundred and eleven verses of different measures, and abounds with fluxional forms and series, to be found in no work of foreign or other Indian countries.'' For more details on the life and work of Sankara Varman one may consult  \cite{Sarma}.

One of \'Sankara Varman's contributions to mathematics was his computation of the value of the mathematical constant $\pi$ correct to 17 decimal places. This in itself, when compared to the computation of $\pi$ by William Shanks correct to 507 decimal places in 1853, is of not much mathematical significance. But, in the scenario and context of the nature of traditional mathematical work in the early nineteenth century India in general and in Kerala in particular, it acquires some significance and indicates that there were persons at that time in Kerala who loved the labour of mathematical computations for its own sake. This paper is an attempt to reproduce in its totality the numerical computations that produced the value of $\pi$ correct to 17 decimal places. 

\section{The Sanskrit verses in {\em Sadratnam\=ala} specifying the value of \texorpdfstring{$\pi$}{pi}}
The verse 2 in Chapter 4 of Sadratnamala specifies the method of computation of the circumference of a circle and also gives the circumference of a circle having diameter $10^{17}$. The verse is reproduced below:
\begin{center}
\includegraphics[width=9cm]{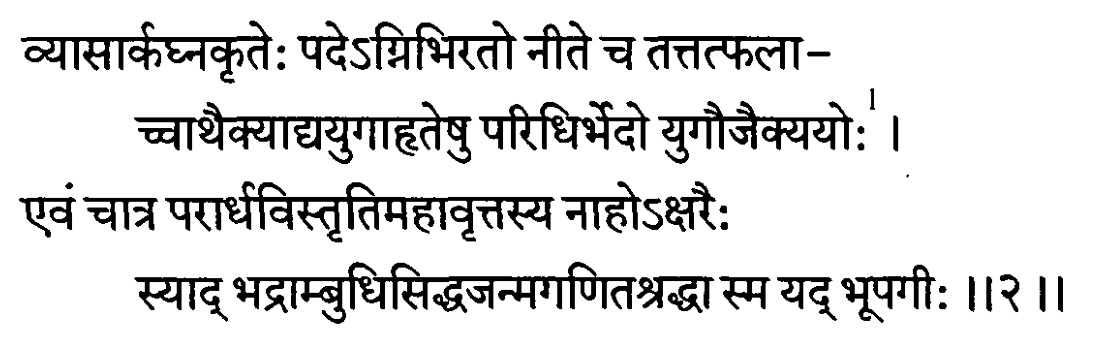}
\end{center}
The verse has been translated into English as follows (see \cite{Madhavan}):
\begin{quote}
\qquad ``Square the diameter of the circle, multiply by 12 and extract the square root. With this as the first term, develop the series thus: Divide this continuously by 3 and divide by 1, 3, 5, 7, 9, 11, ... and form the terms. Add the odd terms and the even terms. Subtract the sum of the even terms from the sum of the odd terms. The result is the circumference of the circle. 

\qquad When this is done, the circumference of the big circle with diameter $10^{17}$ units is 
{\em bhadr\=ambudhisiddha-janma-ga\d{n}ita-\'sradh\=a sma yad bh\=upag\=i\d{h}} 
(by assignment of {\em ka\d{t}apa\=adi} numerals).''
\end{quote}

Let $D$ be the diameter of a circle and $C$ its circumference. Then what the first half of the verse says can be put in the following form:
\begin{align}
C & = \left(\frac{\sqrt{12D^2}}{1} + \frac{\sqrt{12D^2}}{5\cdot 3^2} 
+\frac{\sqrt{12D^2}}{9\cdot 3^4} +\frac{\sqrt{12D^2}}{13\cdot 3^6} + \cdots\right)-\notag \\
&\phantom{=\,\,\,}
\left( \frac{\sqrt{12D^2}}{3\cdot 3} +\frac{\sqrt{12D^2}}{7\cdot 3^3}+\frac{\sqrt{12D^2}}{11\cdot 3^5} + \frac{\sqrt{12 D^2}}{15\cdot 3^7}+\cdots\right)\label{series}
\end{align}
This can be put succinctly in the following form:
\begin{equation}\label{series2}
 C =  \sum_{k=1}^\infty (-1)^{k-1}\frac{\sqrt{12D^2}}{(2k-1)3^{k-1}}
\end{equation}

The second half of the verse gives the value the circumference of a circle of diameter $10^{17}$ ({\em par\=arddha} in Sanskrit). The value is presented in the {\em ka\d{t}apay\=adi} notation.
When the Sanskrit letters in {\em bhadr\=ambudhisiddha-janma-ga\d{n}ita-\'sradh\=a sma yad bh\=upag\=i\d{h}} are replaced by the corresponding numerals, we get the following number:
$$
42\,39\,79\,85\,35\,62\,95\,14\,13
$$
As is well known, the Indian tradition is to write digits from right-to-left which means that the digit in the unit's place is written as the right most digit and the digits with the highest place value is written as the left-most digit. Writing the number in the current mathematical notation we get the number as 
$$
C^\prime = 31\,41\,59\,26\,53\,58\,97\,93\,24.
$$
Sankara Varman claims that is the circumference of a circle whose diameter is one {\em par\=ardha} which is the Sanskrit name for $10^{17}$.

The true value of $\pi$ correct to 20 decimal places is 
$$
\pi = 3.14159265358979323846.
$$
Hence, the circumference of a circle of diameter $10^{17}$ is
\begin{align*}
C & = 3.14159265358979323846\times 10^{17}\\
  & = 31\,41\,59\,26\,53\,58\,97\,93\,23.846 \\
  & \approx  31\,41\,59\,26\,53\,58\,97\,93\,24\\
  & = C^\prime
\end{align*}
Thus the value given by Sankara Varman agrees with the true value of the circumference. 
\subsubsection*{Remark 1}
The result \eqref{series} is not due to \=Sankara Varman. The discovery of the result has been attributed to Sangamagr\=ama Madhava the founder of the Kerala school of astronomy and mathematics (see \cite{MacTutor}, \cite{History} p.24).
\subsubsection*{Remark 2} 
Since $C=\pi D$, from Eq.\eqref{series2}, we have
$$
\frac{\pi}{6} = \sum_{k=1}^\infty (-1)^{k-1}\frac{1}{2k-1}\left(\frac{1}{\sqrt{3}}\right)^{2k-1}.
$$
which is the M\=adhava-Gregory series for $\arctan x$ with $x=\frac{1}{\sqrt{3}}$.
\section{Some issues}
In this context we would like to raise the following questions relating to Sankara Varman's statement in the cited verse.
\subsection{Why a circle of diameter \texorpdfstring{$10^{17}$}{10 followed by 17 zeros}?}
The first question is of only historical curiosity: Why did Sankara Varman consider a circle of diameter $10^{17}$? Is there any special significance to the number $10^{17}$? There should definitely be some significance; he could have easily chosen a circle of diameter, say, $10^{20}$. In several ancient Indian astronomical texts one can see lists of names of powers of ten and these lists stop with the number $10^{17}$ (see \cite{Gupta}, \cite{Datta} p.13). For example, Srídhara (750) gives the following names: {\em eka}, {\em da\'sa}, {em \'sata}, {\em sahasra}, {\em ayuta}, {\em lak\d sa}, {\em prayuta}, {\em ko\d{t}i}, {\em arbuda}, {\em abja}, {\em kharva}, {\em nikharva}, {\em mahâ­saroja}, {\em \'sa\d nku}, {\em sarit\=a-pati}, {\em antya}, {\em madhya}, {\em parârdha}, and adds that the decuple names proceed even beyond this.

Perhaps Sankara Varman chose the number $10^{17}$ because it is the largest named number in Indian astronomical literature. 
\subsection{How did Sankara Varman verify that his value is the true value of the circumference?}
It seems that there is no indication in the available literature that suggests that Sankara Varman did  indeed verify that the value obtained by him is the true value other than perhaps a repeat of the calculations. As we examine below we have no details of how the computations were exactly carried out. Our calculations as described in Section  \label{Computation} show that Sankara Varman obtained the remarkable accuracy of 17 decimal places not by design but by chance.
\subsection{How exactly were the computations performed?}
The computation of the circumference of a circle of diameter $10^{17}$ involves the computation of $\sqrt{12\times (10^{17})^2}$. It may assumed that \'Sankara Varman computed this value using \=Aryabhata's method for the computation of square roots. This method should have been completely familiar to him because {\em \=Aryabhat\=iya} was known and had been thoroughly studied by Kerala astronomers  since at least the time of Sangamagr\=ama M\=adhava (c. 1350 – c. 1425). One of M\=adhava's students, Vata\'{ss}eri Parame\'svaran Namb\=udiri (c. 1380–1460), has written an extensive commentary on {\em \=Aryabhat\=iya}. \=Aryabhata's algorithm for the computation of square roots is also given Jye\d{s}\d{t}hadeva's {\em Yuktibh\=a\d{s}a} (see \cite{Yukti} p.161-162). 

The following verse in {\em \=Aryabhat\=iya} describe the procedure for finding square roots (see \cite{KSShukla} p.36):
\begin{center}
\includegraphics[width=6cm]{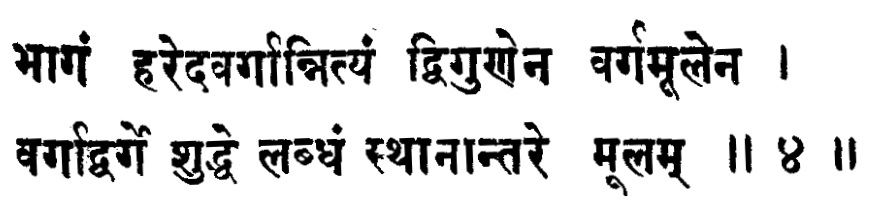}
\end{center}
The verse has been translated as follows:
\begin{quote}
``(Having subtracted the greatest possible square from the last odd place and then having written down the square root of the number subtracted in the line of the square root) always divide the even place (standing on the right) by twice the square root. Then, having subtracted the square (of the quotient) from the odd place (standing on the right), set down the quotient at the next place (i.e., on the right of the number already written in the line of the square root). This is the square root. (Repeat the process if there are still digits on the right.)''
\end{quote}
An example illustrating the application of the above procedure is given in Table 1. From the table we see that 
$$
987654321 = 31426^2 + 60845.
$$
and so $\left\lfloor\sqrt{987654321}\right\rfloor=321426$.
\begin{table}[!t]
\caption{Details of the computation of $\sqrt{987654321}$ by \=Aryabha\d{t}a's method}
\begin{center}
\begin{tabular}{l||l|l}
\hline
{\bf Computations} & {\bf Result} & {\bf Notes}\\
\hline\hline
{\bf 9,8 7,6 5,4 3,2 1} & & \\
\hline\hline
{\bf 9} $-$ &{\bf 3} &$\lfloor \sqrt 9 \rfloor = \mathbf 3$ \\
9 & & $\mathbf 3^2=1$\\
\hline
0 {\bf 8} $-$&3 {\bf 1} & $\lfloor 8/(2\cdot 3)\rfloor = \mathbf 1$\\
0 6 &    & $\mathbf 1\cdot(2\cdot 3)$\\
\hline
\phantom{0 }2 {\bf 7} $-$&  & \\
\phantom{2 6 }1 &  & $\mathbf 1^2=1$\\
\hline
\phantom{0 }2 6 {\bf 6} $-$ &3 1 {\bf 4}     & $\lfloor 266/(2\cdot 31)\rfloor =\mathbf  4$\\
\phantom{0 }2 4 8 &     & $\mathbf 4\cdot(2\cdot 31) = 248$
\\
\hline
\phantom{0 0 }1 8 {\bf 5} $-$ & & \\
\phantom{0 0 0 }1 6  & & $\mathbf 4^2=16$\\
\hline
\phantom{0 0 }1 6 9 {\bf 4} $-$ &3 1 4 {\bf 2}  & $\lfloor 1694/(2\cdot 314)\rfloor = \mathbf 2 $\\
\phantom{0 0 }1 2 5 6& & $\mathbf 2\cdot(2\cdot 314)=1256$\\
\hline
\phantom{0 0 0 }4 3 8 {\bf 3} $-$&  & \\
\phantom{0 0 0 0 0 0 }4 & & $\mathbf 2^2=4$\\
\hline
\phantom{0 0 0 }4 3 7 9 {\bf 2} $-$ & 3 1 4 2 {\bf 6} &$\lfloor 43792/(2\cdot 3142)\rfloor = \mathbf 6$\\
\phantom{0 0 0 }3 7 4 0 4 & & $\mathbf 6\cdot(2\cdot 3142) = 37704$\\
\hline
\phantom{0 0 0 0 0 }6 0 8 8 {\bf 1} &  & \\
\phantom{0 0 0 0 0 0 0 0 }3 6 & & $\mathbf 6^2=36$\\
\hline
\phantom{0 0 0 0 0 }6 0 8 4 5    &  & Remainder\\
\hline
\end{tabular}
\end{center}
\end{table}

\subsection{Handling fractional parts}
The computations involve, besides the operations of addition, subtraction and multiplication, the operation of extraction of square root and division. The square root to be extracted is $\sqrt{12\times (10^{17})^2}$. The number $12\times (10^{17})^2$ is not an exact square and hence its square root is not an exact integer. So some approximate value has to be used for further numerical computations. The same is the case with divisions. The result of the division of an integer by another integer is more often not an exact integer. We can carry out the computations in several possible ways.
\begin{enumerate}
\item
{\bf Ignore the fractional part at every operation}: We may perform the computation by ignoring the fractions at every stage of the computation. This means that that the result of an arithmetic operation involving positive integers shall always be taken as the integral part of the result of the operation. Thus, in this case we shall assume that the value of $\sqrt{12\times (10^{17})^2}$ is $\lfloor \sqrt{12\times (10^{17})^2}\rfloor$ where $\lfloor\,\,\,\rfloor$ denotes the floor function.

Using the floor function, this computational scheme can be specified as follows:
\begin{align*}
x_1 & = \left \lfloor \sqrt{12D^2} \right \rfloor\\
x_k & = \left\lfloor \frac{x_{k-1}}{3}\right\rfloor\text{ for } k=2,3,4, \cdots\\
t_k & = \left\lfloor \frac{x_k}{2k-1}\right\rfloor \text{ for } k =1, 2, ,3 ,4 , \cdots\\
O & = t_1+t_3+t_5+\cdots\\
E & = t_2+t_4+t_6+\cdots\\
C & = O-E
\end{align*}
\item
{\bf Round the result of an operation to the nearest integer}: Since, currently, there are several different methods for rounding of numbers we may consider one of the most commonly used methods of rounding which is rounding to the nearest integer. If $x$ is a positive real number, the result of rounding $x$ to the nearest integer can be expressed as $\lfloor x+\frac{1}{2} \rfloor$. One may round off the result of each division to nearest integer and take the rounded value as the result of division.

Using the floor function, this computational scheme can be specified as follows:
\begin{align*}
x_1 & = \left\lfloor \sqrt{12D^2} +\frac{1}{2}\right\rfloor\\
x_k & = \left\lfloor \frac{x_{k-1}}{3} +\frac{1}{2}\right\rfloor\text{ for } k=2,3,4, \cdots\\
t_k & = \left\lfloor \frac{x_k}{2k-1} +\frac{1}{2} \right\rfloor \text{ for } k =1, 2, ,3 ,4 , \cdots\\
O & = t_1+t_3+t_5+\cdots\\
E & = t_2+t_4+t_6+\cdots\\
C & = O-E
\end{align*}
\item
{\bf Do the rounding at the end}: In this method, one may retain the fractional parts in all the terms and do the rounding at the last computation. However, to start with, we have to assign a value for $\sqrt{12\times (10^{17})^2}$.

\begin{alignat*}{2}
x_1 & = \sqrt{12D^2} \\
x_k & = \frac{x_{k-1}}{3}&& \text{ for } k=2,3,4, \cdots\\
t_k & = \frac{x_k}{2k-1}&&  \text{ for } k =1, 2, ,3 ,4 , \cdots\\
O & = t_1+t_3+t_5+\cdots&&\\
E & = t_2+t_4+t_6+\cdots&&\\
C & = \lfloor O-E \rfloor &&\text{ignoring fractional part}\\
\text{or}\qquad C & =\left\lfloor O-E + \tfrac{1}{2}\right\rfloor&&\text{rounding to nearest integer}
\end{alignat*}
Even though this is a theoretical possibility, it is difficult to believe that \'Sankara Varman did his computation using this method of rounding.
\end{enumerate}
\section{Computations}
\subsection{Ignore the fractional part at every operation}

\begin{enumerate}[{\bf Step }\bf 1. ]
\item
Diameter of the circle is $D=10^{17}$.
\item
Square the diameter of the circle to get $D^2=(10^{17})^2=10^{34}$.
\item
Multiply by $12$ to get $12D^2=12\times 10^{34}$.
\item
Find the square root of $12D^2$. Applying \=Aryabhata's algorithm and ignoring the fractional part we get 
$$
\sqrt{12D^2} = 346410161513775458
$$
\item 
With this as the first term, develop the series thus: Divide this continuously by 3 and divide by 1, 3, 5, 7, 9, 11, ... and form the terms. The results of continuous division by $3$ are shown the second column of Table 1 and the results of continuous division by 1, 3, 5, ... are shown in the last column of the same table. We have computed $38$ terms of the series in Eq.\eqref{series2} because the operations of successive divisions by $3$ produce the value $0$ for the $38$-th term as can be seen from Table 1 and further divisions by $3$ will only give 0's.
\item
The following sums are calculated.
\begin{align*}
\text{Sum of odd numbered terms } O & = 346410161513775458 + \\
&\phantom{ =0 } 1832857997427383 + \\
&\phantom{ =0 } 129596020020118 + \cdots\\
& = 354623317218212158 \\
\text{Sum of even numbered terms } E & = 7698003589195010 + \\
&\phantom{ =0 } 475185406740432 + \\
&\phantom{ =0 } 36552723595417 +\cdots\\
& = 40464051859232834 
\end{align*}

\item
The circumference $C$ of a circle of diameter $10^{17}$ is given by
\begin{align*}
C & = O - E\\
& = 354623317218212158 - 40464051859232834 \\
& = 314159265358979324
\end{align*}
This is exactly the value given by Sankara Varman.
\end{enumerate}
\subsubsection*{Remark}
All the above computations can be done by hand and \'Sankara Varman would have performed the computations  definitely by hand and these computations would not have been very cumbersome and tedious. 
However, for a modern scholar or a student doing these computations would surely appear extremely complicated and highly prone to error. 
Now they can verify the accuracy of the computations by using any one of the several computer algebra systems available in the market because such systems are good at performing exact arithmetic operation with arbitrarily large integers.  
The author of the present paper did actually verify the results using the {\em Maxima} (version 5.46.0) package which is a free software released under the terms of the GNU General Public License. The {\em Maxima}  code for reproducing the results is given in Appendix 2.
\smallskip
\newpage
%
%\begin{table}
{Table 1: Computation of the terms in the series for the circumference}\quad\\[2mm]
\quad\hspace{-0.1in}
\begin{tabular}{c|r|c|c|l}
\hline
{Term}&{After dividing} &{Division}& {Sign} & {After division}\\
{position}  & {previous term by 3} &{by} ${2k-1}$&  & {by} ${2k-1}$\\
($k$)  & ($x_k$) &  &  & ($t_k$)\\
\hline\hline
1  &346410161513775458&$(\div 1) \rightarrow$& $+$ & 346410161513775458\\
2&115470053837925152 &$(\div 3) \rightarrow$& $-$ & 38490017945975050\\
3&38490017945975050 &$(\div 5) \rightarrow$& $+$ & 7698003589195010\\
4 & 12830005981991683 &$(\div 7) \rightarrow$& $-$ & 1832857997427383\\
5 & 4276668660663894 &$(\div 9) \rightarrow$& $+$ & 475185406740432\\
6 & 1425556220221298 &$(\div 11) \rightarrow$& $-$ & 129596020020118\\
7 & 475185406740432 &$(\div 13) \rightarrow$& $+$ & 36552723595417\\
8 & 158395135580144 &$(\div 15) \rightarrow$& $-$ & 10559675705342\\
9 & 52798378526714 &$(\div 17) \rightarrow$& $+$ & 3105786972159 \\
10 & 17599459508904 &$(\div 19) \rightarrow$& $-$ & 926287342573\\
11 & 5866486502968 &$(\div 21) \rightarrow$& $+$ & 279356500141\\
12 & 1955495500989 &$(\div 23) \rightarrow$& $-$ & 85021543521\\
13 & 651831833663 &$(\div 25) \rightarrow$& $+$ & 26073273346\\
14 & 217277277887 &$(\div 27) \rightarrow$& $-$ & 8047306588\\
15 & 72425759295 &$(\div 29) \rightarrow$& $+$ & 2497439975\\
16 & 24141919765 &$(\div 31) \rightarrow$& $-$ & 778771605\\
17 & 8047306588 &$(\div 33) \rightarrow$ & $+$ & 243857775\\
18 & 2682435529 &$(\div 35) \rightarrow$& $-$ & 76641015\\
19 & 894145176 &$(\div 37) \rightarrow$ & $+$ & 24166085\\
20 & 298048392 &$(\div 39) \rightarrow$& $-$ & 7642266\\
21 & 99349464 &$(\div 41) \rightarrow$ & $+$ & 2423157\\
22 & 33116488 &$(\div 43) \rightarrow$ & $-$ & 770150\\
23 & 11038829 &$(\div 45) \rightarrow$ & $+$ & 245307\\
24 & 3679609 &$(\div 47) \rightarrow$ & $-$ & 78289\\
25 & 1226536 &$(\div 49) \rightarrow$ & $+$ & 25031\\
26 & 408845 &$(\div 51) \rightarrow$ & $-$ & 8016\\
27 & 136281&$(\div 53) \rightarrow$ & $+$ & 2571\\
28 & 45427 &$(\div 55) \rightarrow$ & $-$ & 825\\
29 & 15142 &$(\div 57) \rightarrow$ & $+$ & 265\\
30 & 5047 &$(\div 59) \rightarrow$& $-$ & 85\\
31 & 1682 &$(\div 61) \rightarrow$& $+$ & 27\\
32 & 560 &$(\div 63) \rightarrow$& $-$ & 8\\
33 & 186 &$(\div 65) \rightarrow$ & $+$ & 2\\
34 & 62 &$(\div 67) \rightarrow$ & $-$ & 0\\
35 & 20 &$(\div 69) \rightarrow$ & $+$ & 0\\
36 & 6 &$(\div 71 ) \rightarrow$& $-$ & 0\\
37 & 2 &$(\div 73) \rightarrow$ & $+$ & 0\\
38 & 0 &$(\div 75) \rightarrow$& $-$ & 0\\
\hline
\end{tabular}
%\end{table}
%\end{center}
%
\subsection{Rounding to the nearest integer at every stage of computation}
In this method computation we have (this may be verified by making appropriate changes in the {\em Maxima} code given in Appendix 2)
\begin{align*}
\sqrt{12D^2} & = 346410161513775459\\
O & = 354623317218212169 \\
E & = 40464051859232844 \\
C & = O - E \\
  & = 314159265358979325
\end{align*}
It may be noted that this value is different from the value given by Sankara Varman. It differs by 1 in the unit's place.
\subsection{Do the rounding at the end}
In this case (again, this may be verified using {\em Maxima} software) we have
\begin{alignat*}{2}
C & =314159265358979323.84375&&\\
  & \approx 314159265358979323 &&\quad\text{(ignoring fractional part)}\\
  & \approx 314159265358979324&&\quad \text{(rounding to nearest integer)}
  \end{alignat*}
The above result was also obtained by considering the first $38$ terms of the series Eq.\eqref{series2}.
\section{Sangamagr\=ama M\=adhava}
Sangamagr\=ama M\=adhava, referred to as M\=adhavan Ila\~{nn}ppa\d{ll}i Empr\=an in a certain manuscript, was a member of the E\d{m}pr\=antiri community which is a community of Tu\d{l}u speaking Br\=ahmins settled in Kerala.  His forefathers are believed to have migrated from the Tulu land to settle in Kudallur village, which is situated on the southern banks of the Nila river not far from Tirunnavaya, a generation or two before his birth and lived in a house known as Ila\~{nn}ppa\d{ll}i whose present identity is unknown (see \cite{Divakaran}). There are no definite evidences to pinpoint the period during which M\=adhava flourished. From circumstantial evidences historians have assigned the date c. 1340 - 1425 to M\=adhava. 
\section{M\=adhava's value of \texorpdfstring{$\pi$}{pi}}

A verse ascribed to M\=adhava and quoted in the {\em Kriyakramakri} commentary of {\em L\=il\=avati} by \'Sankara specifies  a numerical value for $\pi$. The numerical value is encoded in the {\em bh\=utasa\d{m}khy\=a} system. In this system the numerals are expressed by names of things, beings or concepts, which, naturally or in accordance with the {\em \'S\=astras}, connote numbers; for example, the number two may be denoted by anything which suggests the number 2. (For a list of words commonly used to denote numbers see \cite{Datta} p.14-17.)   

The verse and an English translation are reproduced below:

\begin{center}
\includegraphics[width=8cm]{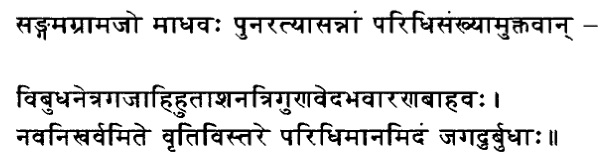}
\end{center}
\begin{quote}
``The teacher M\=adhava also mentioned a value of the circumference
closer [to the true value] than that:

{\em vibudha} [gods = 33] {\em netra} [eyes = 2] {\em gaja} [elephants = 8] {\em ahi} [serpents = 8] {\em hut\=a\'sana} [fires = 3] {\em tri} [three = 3]
{\em guna} [qualities = 3] {\em veda} [Vedas = 4] {\em bha} [{\em nak\d satras} = 27] {\em v\=arana} [elephants = 8] {\em b\=ahav\=a\d{h}} [arms = 2 is the
measure of a circle of diameter {\em nava} [9] {\em nikharva} [100,000,000,000].''

\end{quote}

Thus, the circumference of a circle of diameter $9\times 10^{11}$ in the Indian notational scheme is 
$$
33\,\,2\,\,8\,\,8\,\,3\,\,\,3\,\,4\,\,27\,\,8\,\,2
$$
and in modern notational scheme is 
$$
2,827,433,388,233.
$$
With this value for the circumference we have 
\begin{align}
\pi &\approx  \frac{2,827,433,388,233}{9\times 10^{11}}\notag\\
   &\approx \underline{3.1415926535}\,\,\underline{92}22222\ldots\label{pi1}
\end{align}
and the true value of $\pi$ is
$$
\pi = 3.1415926535\,\,89793238\ldots
$$
The approximate value of $\pi$ is correct up to the $10$-th decimal place. 

With the true value of $\pi$, the circumference of a circle of diameter $9\times 10^{11}$ is
\begin{align*}
\pi \times 9\times 10^{11} 
& = 2827433388230.81391461637904495\ldots\\
& = 2827433388231 \qquad \text{rounding to nearest integer}\\
& = 2827433388230  \qquad \text{ignoring the fractional part}
\end{align*}
With these approximate values for the circumference of a circle of diameter $9\times 10^{11}$ we have
\begin{align}
\pi &\approx \frac{2827433388231}{9\times 10^{11}} = \underline{3.1415926535}\,\,\underline{90} \qquad\text{(this is exact)}\\
\pi &\approx \frac{2827433388230}{9\times 10^{11}} =
\underline{3.1415926535}\,\,\underline{88}88888\ldots
\end{align}
Both these values are closer to the true value of $\pi$ than the approximate value given in Eq.\eqref{pi1}.
\section{How M\=adhava might have computed the circumference?}
\'Sankara Varman has indicated that he computed the value of $\pi$ using the series in Eq.\eqref{series}. But neither N\-ilakan\d{t}ha Somay\=aji nor \'Sankara V\=ariar has indicated how exactly M\=adhava calculated the circumference of a circle of diameter $9\times 10^{11}$. In this section we examine some of the possible methods M\=adhava might have employed to calculate the circumference. 

The following formulae for computing the circumference of a circle of diameter $D$ have been attributed to M\=adhava.
\begin{enumerate}[{\bf F}\bf 1.]
\item
$
\left(\dfrac{\sqrt{12D^2}}{1} + \dfrac{\sqrt{12D^2}}{5\cdot 3^2} 
+\dfrac{\sqrt{12D^2}}{9\cdot 3^4}  + \cdots\right)-
\left( \dfrac{\sqrt{12D^2}}{3\cdot 3} +\dfrac{\sqrt{12D^2}}{7\cdot 3^3}+\dfrac{\sqrt{12D^2}}{11\cdot 3^5} +\cdots\right)\label{seriesX}
$
This is same as the series in Eq.\eqref{series} and it is the series used by \'Sankara Varman to compute his value of $\pi$.
\item
$\dfrac{4D}{1}-\dfrac{4D}{3}+\cdots + (-1)^{n-1}\dfrac{4D}{2n-1}+(-1)^{n}4D F(n)$

Here $F(n)$ is a correction term used to give better approximation to the circumference. Three expressions for $F(n)$ have been ascribed to M\=adhava the third being the one which gives the best approximation to the circumference.
\begin{enumerate}
\item
$F_1(n)=\dfrac{1}{4n}$
\item
$F_2(n) = \dfrac{n}{4n^2+1}$
\item
$F_3(n)= \dfrac{n^2+1}{n(4n^2+5)}$
\end{enumerate}
\item
$ 3D + \dfrac{4D}{3^3-3}-\dfrac{4D}{5^3-5}+\dfrac{4D}{7^3-7} - \cdots$
\item
$\dfrac{16D}{1^5 + 4\cdot 1} - \dfrac{16D}{3^5+4\cdot 3}+\dfrac{16D}{5^5+4\cdot 5} - \dfrac{16D}{7^5+4\cdot 7}+ \cdots$
\end{enumerate}
\subsection{Results of computations}
\subsubsection{Formula F1}
\begin{center}
Table 2: Circumference of a circle of diameter $9\times 10^{11}$ using formula {\bf F1}\\[2mm]
\begin{tabular}{c||c|c|c}
\hline
      &\multicolumn{3}{c}{Circumference ($C$)} \\
      \cline{2-4}
Number of & Ignoring & Operation-wise & Final \\
terms ($n$) & fractions & rounding & rounding\\
\hline\hline
18     &       2827433388065     &       2827433388065            &2827433388066   \\
19     &       2827433388283     &       2827433388282            &2827433388283   \\
20       &     2827433388214     &       2827433388214            & 2827433388214   \\
21&            2827433388236    &        2827433388235            &2827433388236   \\
22 &           2827433388229      &      2827433388229            &2827433388229   \\
23  &          2827433388231      &      \underline{2827433388231} &2827433388231 \\  
24   &        \underline{2827433388230} &       2827433388231            &\underline{2827433388230} \\ 
25    &        2827433388230    &        2827433388231            &2827433388230  \\ 
26     &       2827433388230     &       2827433388231            &2827433388230 \\  
27      &      2827433388230      &      2827433388231            &2827433388230 \\  
\hline
\end{tabular}
\end{center}

Under the specified conditions, the sum of the infinite series in {\bf F1} is the underlined number in Table 2. Note that the series does not sum to M\=adhava's value for the circumference; instead, it sums to a value much closer to the true value of the circumference.
\subsubsection{Formula F2}
\begin{center}
Table 3: Circumference of a circle of diameter $9\times 10^{11}$ using formula {\bf F2} with correction term $F_3(n)$\\[2mm]
\begin{tabular}{c||c|c|c}
\hline
  & \multicolumn{3}{c}{Circumference ($C$)}\\
  \cline{2-4}
Number of & Ignoring & Operation-wise & Final \\
terms ($n$)& fractions & rounding & rounding\\
\hline\hline
35        &    2827433388224          & 2827433388221            &2827433388223   \\
36         &   2827433388238&           2827433388235            &2827433388237   \\
37          &  2827433388227 &          2827433388223            &2827433388225   \\
38           & 2827433388235  &         {\bf 2827433388233}            &2827433388235   \\
39&            2827433388228   &        2827433388225            &2827433388227   \\
40 &           2827433388234    &       2827433388231            &{\bf 2827433388233}   \\
41  &          2827433388229     &      2827433388225            &2827433388228   \\
42   &         2827433388234      &     2827433388230            &2827433388232   \\
43    &        2827433388230       &    2827433388225            &2827433388228   \\
44     &       {\bf 2827433388233}        &   2827433388229            &2827433388232   \\
45      &      2827433388231         &  2827433388226            &2827433388229   \\
46       &     {\bf 2827433388233}          & 2827433388229            &2827433388231   \\
47        &    2827433388232           &2827433388226            &2827433388229   \\
48         &   2827433388232&           2827433388228            &2827433388231   \\
49          &  2827433388231 &          2827433388227            &2827433388230   \\
50           & 2827433388232  &         2827433388228            &2827433388231   \\
51&            2827433388231&           2827433388226            &2827433388230 \\
52 &           2827433388232 &          2827433388227            &2827433388231   \\
53    &        2827433388232   &        2827433388226            &2827433388230   \\
54  &          2827433388231   &        2827433388227            &2827433388231   \\
55   &         2827433388231    &       2827433388226            &2827433388230   \\
56    &        2827433388231     &      2827433388227            &2827433388231   \\
57     &       2827433388231      &     2827433388227            & \underline{2827433388230}   \\
58      &      2827433388230       &    2827433388228            &2827433388231   \\
59       &     2827433388231        &   2827433388227            &2827433388230   \\
60        &    2827433388230         &  2827433388228            &2827433388230   \\
61         &   2827433388230          & 2827433388228            &2827433388230   \\
62          &  2827433388231           &2827433388228            &2827433388230   \\
63           & 2827433388231&           2827433388228            &2827433388230   \\
64            &2827433388231  &         2827433388228            &2827433388230   \\
65&            2827433388232   &        2827433388228            &2827433388230   \\
\hline
\end{tabular}
\end{center}

The numbers in boldface indicate the positions where M\=adhava's value for the circumference appear. But, under any of the three methods of rounding, the series does to sum to M\=adhava's value. Also, if the computations are done with the rounding at the final stage then the sum converges to the value $2827433388230$ (the underlined number in the table) which is different from M\=adhava's value but closer to the true value than M\=adhava's value. The other rounding procedures do not produce a repeating value even if we do computations with number of terms up to a million. Moreover, if we take a million terms, the term-by-rounding yields the sum $2827433387851$ which is a very poor approximation to the true value of the circumference. 
\subsubsection{Formula F3}
If we apply this formula, then the series indeed sum to a definite value. Computations using {\em Maxima} software give the following results:
\begin{itemize}
\item
If we ignore all fractions then $C=2827433388211$ for all $n\ge 7663$. Thus in this case the series sums to $2827433388211$. The number $7663$ is the smallest integer $n$ for which 
$$
\frac{4\times\times 9\times 10^{11}}{(2n-1)^3 -(2n-1)} \ge 1.
$$
\item
If we apply rounding at each operation then $C=2827433388236$ for $n\ge 9655$. SO, in this case the series sums to $2827433388236$. The number $9655$ is the smallest integer $n$ for which 
$$
\frac{4\times\times 9\times 10^{11}}{(2n-1)^3 -(2n-1)} + \frac{1}{2}\ge 1.
$$
\item
If the rounding is done only at the last stage then
$C = 2827433388231$ for $n\ge 8950$ and hence in this case the series sums to $2827433388231$.
\end{itemize}

Thus formula {\bf F3} also  does not produce M\=adhava's value for the circumference.
\subsubsection{Formula F4}

\begin{center}
Circumference of a circle of diameter $9\times 10^{11}$ using formula {\bf F4}\\[2mm]
\begin{tabular}{c||c|c|c}
\hline
 & \multicolumn{3}{c}{Circumference ($C$)}\\ \cline{2-4}
Number of & Ignoring & Operation-wise & Final \\
terms ($n$)& fractions & rounding & rounding\\
\hline\hline
210&            2827433388226&           2827433388233            &2827433388230 \\  
211&            2827433388227&           2827433388234            &2827433388231   \\
212&            2827433388226 &          2827433388233            &2827433388230   \\
213&            2827433388227   &        2827433388234            &2827433388231   \\
214&            \underline{2827433388226}   &        2827433388233            &2827433388230   \\
215&            2827433388226    &       2827433388234            &2827433388231   \\
216&            2827433388226     &      2827433388233            &2827433388230   \\
217&            2827433388226      &     2827433388234            &2827433388231   \\
218&            2827433388226       &    2827433388233            &2827433388230   \\
219&            2827433388226        &   2827433388234            &2827433388231   \\
220&            2827433388226         &  2827433388233            &2827433388230   \\
221&            2827433388226          & 2827433388234            &2827433388231   \\
222&            2827433388226           &2827433388233            &2827433388230   \\
223&            2827433388226&           2827433388234            &2827433388231   \\
224&            2827433388226 &          2827433388233            &2827433388230   \\
225&            2827433388226  &         2827433388234            &2827433388231   \\
226&            2827433388226   &        2827433388233            &2827433388230   \\
227&            2827433388226    &       2827433388234            &2827433388231   \\
228&            2827433388226     &      2827433388233            &2827433388230   \\
229&            2827433388226      &     2827433388234            &2827433388231   \\
230&            2827433388226       &    2827433388233            &2827433388230   \\
231 &           2827433388226        &   2827433388234            &2827433388231   \\
232  &          2827433388226         &  2827433388233            &2827433388230   \\
233   &         2827433388226          & 2827433388234            &2827433388231   \\
234    &        2827433388226           &2827433388233            &2827433388230   \\
235     &       2827433388226&           2827433388234            &\underline{2827433388231}   \\
236      &      2827433388226 &          2827433388233            &2827433388231   \\
237       &     2827433388226  &         2827433388234            &2827433388231   \\
238      &      2827433388226   &        2827433388233            &2827433388231   \\
239       &     2827433388226    &       2827433388234            &2827433388231   \\
240        &    2827433388226     &      2827433388233            &2827433388231   \\
241         &   2827433388226      &     2827433388234            &2827433388231   \\
242          &  2827433388226       &    2827433388233            &2827433388231   \\
243           & 2827433388226        &   2827433388234            &2827433388231   \\
244            &2827433388226         &  2827433388233            &2827433388231   \\
245&            2827433388226          & 2827433388234            &2827433388231   \\
246 &           2827433388226           &\underline{\bf 2827433388233}            &2827433388231   \\
247  &          2827433388226&           2827433388233            &2827433388231   \\
248   &         2827433388226 &          2827433388233            &2827433388231   \\
249    &        2827433388226  &         2827433388233            &2827433388231   \\
250     &       2827433388226   &        2827433388233            &2827433388231\\
\hline
\end{tabular}
\end{center}

As in the case of formula {\bf F3}, we note the following:
\begin{itemize}
\item
If we ignore the fractions the series in {\bf F4} sums to $2827433388226$ which is different from M\=adhava's value.
\item
If we perform rounding for each computation, the series in {\bf F4} sums to $2827433388233$ which is precisely M\=adhava's value for the circumference.
\item
If we perform rounding at the final stage only, the series in {\bf F4} sums to $2827433388231$ which is different from M\=adhava's value but closer to the true value of the conference.
\end{itemize}
\subsection{Conclusion: M\=adhava erred in his computations?}
From the numerical computations described in Section 6.1 we arrive at the following conclusions.
\begin{itemize}
\item
Among the four formulas considered only the following formula has produced M\=adhava's value for the circumference:
$$
\dfrac{16D}{1^5 + 4\cdot 1} - \dfrac{16D}{3^5+4\cdot 3}+\dfrac{16D}{5^5+4\cdot 5} - \dfrac{16D}{7^5+4\cdot 7}+ \cdots
$$
\item
This series produces M\=adhava's value for the circumference only if the computations are performed with rounding to the nearest integer at every stage of the computation.
\item
To ensure that the series indeed yields M\=adhava's value, the summation is to be performed at leat up to the $247$-th term.
\end{itemize}
Rather than concluding that M\=adhava arrived at his value of $2827433388233$ for the circumference as described above, a much simpler explanation might be that M\=adhava erred in his computation of the circumference of a circle of radius $9\times 10^{11}$ however disquieting this explanation might be!

%%%%%%%%%%%%%%%%%%%%%%%%

%
\end{document}